\documentclass{commat}

\usepackage{mathtools}

\DeclarePairedDelimiter\ffact{\langle}{\rangle}
\DeclarePairedDelimiter\ffactb{\Big\langle}{\Big\rangle}

\newcommand{\genstirlingII}[3]{%
  \genfrac{\{}{\}}{0pt}{#1}{#2}{#3}%
}
\newcommand{\stirlingsk}[2]{\genstirlingII{}{#1}{#2}}

\newcommand{\genstirlingI}[3]{%
  \genfrac{[}{]}{0pt}{#1}{#2}{#3}%
}
\newcommand{\stirlingfk}[2]{\genstirlingI{}{#1}{#2}}

\newcommand{\C}{\mathbb{C}}

\title{%
    Solution of the equation $y^{\prime}=f(y)$ and Bell Polynomials
    }

\author{%
    Ronald Orozco L\'opez
    }

\affiliation{
    \address{Department of Mathematics 
    Universidad de los Andes, 111711, Bogot\'a Colombia 
        }
    \email{%
    rj.orozco@uniandes.edu.co
    }
    }

\abstract{%
    In this paper we use Fa\`a di Bruno's formula to associate Bell polynomial values to differential equations of the form $y^{\prime}=f(y)$. That is, we use partial Bell polynomials to represent the solution of such an equation and use the solution to compute special values of partial Bell polynomials.
    }



\VOLUME{31}
\YEAR{2023}
\NUMBER{1}
\firstpage{103}
\DOI{https://doi.org/10.46298/cm.10278}

\begin{paper}

\section{Introduction}

It is a known fact that Bell polynomials are closely related to the derivatives of the composition of functions. For example, Fa\`a di Bruno \cite{Faa}, Foissy \cite{Foissy}, and Riordan \cite{Riordan_j} proved that Bell polynomials are a very useful tool in mathematics to represent the $n$-th derivative of the composition of functions. Also, Bernardini and Ricci \cite{Bernardini}, Yildiz et al. \cite{Yildiz}, Caley \cite{Caley}, and Wang \cite{Wang} showed the relationship between Bell polynomials and differential equations. On the other hand, Orozco \cite{Orozco} studied the convergence of the analytic solution of the autonomous differential equation $y^{(k)}=f(y)$ using Fa\`a di Bruno's formula. We can then consider differential equations as a source for researching special values of Bell polynomials.

In this paper we consider the solution $y(t,x)=g(t+g^{-1}(x))$ of the autonomous differential equation $y^{\prime}=f(y)$ and show how to express this solution by means of Bell polynomials. This will then be used to find special values of partial Bell polynomials. Here we will not consider convergence issues, but formal solutions of such a differential equation.

This paper is organized as follows. We start with basic results on partial and complete Bell polynomials and special values of these. In the third section we show what condition $g(x)$ must satisfy for $g(t+g^{-1}(x))$ to be a solution of $y^{\prime}=f(y)$. We conclude by showing the relationship between Bell polynomials and the solution of the differential equation $y^{\prime}=f(y)$ when $f(x)=ax$, $f(x)=e^{x}$, $f(x)=\sqrt{1-x^2}$, $f(x)=1+x^{2}$, $f(x)=\sqrt{x^{2}\pm1}$, $f(x)=1-x^2$ and $f(x)=\frac{1}{(x-1)^{\alpha-1}}$, where $a,\alpha\in\C$.

\section{Preliminaries}
The following basic results can be found in Comtet \cite{Comtent}, and Riordan \cite{Riordan_b}. Exponential Bell polynomials are used to encode information about the ways in which a set can be partitioned, making them a very useful tool in combinatorial analysis. Bell polynomials are obtained from the derivatives of composite functions and are given by the formula of Fa\`a Di Bruno \cite{Faa}. Bell \cite{Bell}, Gould \cite{Gould_Q}, Mihoubi \cite{Mihoubi}, Wang \cite{Wang2} and Feng Qi \cite{FengQi1}, \cite{FengQi2}, \cite{FengQi3} (among many others) provided important results on these polynomials. We start with the definition of the partial Bell polynomials.

\begin{definition}
The exponential partial Bell polynomials are the polynomials
\begin{equation*}
    B_{n,k}(x_{1},x_{n},\ldots,x_{n-k+1})
\end{equation*}
in the infinite variables $x_{1},x_{2},\ldots$ defined by the series expansion
\begin{equation}\label{eqn_egf}
\exp\left(u\sum_{j=1}^{\infty}x_{j}\frac{t^{j}}{j!}\right)=1+\sum_{n=1}^{\infty}\frac{t^{n}}{n!}\sum_{k=1}^{n}u^{k}B_{n,k}(x_{1},x_{2},\ldots,x_{n-k+1}).
\end{equation}
or equivalently defined by the series expansion of the $k$-th power
\begin{equation}\label{eqn_egf_pow}
    \frac{1}{k!}\left(\sum_{j=1}^{\infty}x_{j}\frac{t^j}{j!}\right)^{k}=\sum_{n=k}^{\infty}B_{n,k}(x_{1},x_{2},\ldots,x_{n-k+1})\frac{t^n}{n!}.
\end{equation}
\end{definition}

The following result gives the explicit way to calculate the partial Bell polynomials
\begin{theorem}
The partial or incomplete exponential Bell polynomials are given by
\begin{equation*}
B_{n,k}(x_{1},\ldots,x_{n-k+1})=\sum\frac{n!}{c_{1}!c_{2}!\cdots c_{n-k+1}!}\left(\frac{x_{1}}{1!}\right)^{c_{1}}\cdots\left(\frac{x_{n-k+1}}{(n-k+1)!}\right)^{c_{n-k+1}},
\end{equation*}
where the summation takes place over all integers $c_{1},c_{2},\ldots,c_{n-k+1}\geq0$, such that
\begin{eqnarray*}
c_{1}+2c_{2}+\cdots+(n-k+1)c_{n-k+1}&=&n,\\
c_{1}+c_{2}+\cdots+c_{n-k+1}&=&k.
\end{eqnarray*}
\end{theorem}

Some values of partial Bell polynomials are
\begin{align*}
    B_{n,k}(0!,1!,\ldots,(n-k)!)&=\stirlingfk{n}{k}\ \ (\textit{Unsigned Stirling number of first kind}),\\
    B_{n,k}(1!,\ldots,(n-k)!)&=\binom{n-1}{k-1}\frac{n!}{k!}\ \ (\textit{Lah number}),\\
    B_{n,k}(1,1,\ldots,1)&=\stirlingsk{n}{k}\ \ (\textit{Stirling number of second kind}),\\
    B_{n,k}(1,2,\ldots,n-k+1)&=\binom{n}{k}k^{n-k}\ \ (\textit{Idempotent number}).
\end{align*}
Then we can see the beautiful relationship that exists between Bell polynomials and numbers like the above.

Feng Qi \cite{FengQi3} deduced the following identity that will be very useful to us
\begin{equation*}\label{eqn_1x_1}
    B_{n,k}(x,1,0,0,\ldots,0)=\frac{1}{2^{n-k}}\frac{n!}{k!}\binom{k}{n-k}x^{2k-n},
\end{equation*}
that together with the identity
\begin{equation}\label{eqn_hom}
    B_{n,k}(abx_{1},ab^{2}x_{2},\ldots,ab^{n-k+1}x_{n-k+1})=a^{k}b^{n}B_{n,k}(x_{1},x_{2},\ldots,x_{n-k+1})
\end{equation}
leads us to
\begin{equation}\label{eqn_ax_a}
    B_{n,k}(ax,a,0,0,\ldots,0)=\frac{a^{k}}{2^{n-k}}\frac{n!}{k!}\binom{k}{n-k}x^{2k-n}.
\end{equation}

Finally we show Fa\`a di Bruno's formula. Let $f$ and $g$ be functions with exponential generating functions $\sum a_{n}\frac{x^{n}}{n!}$ and $\sum b_{n}\frac{x^{n}}{n!}$ 
respectively, with $a_{n},b_{n}\in \C$. Then
\begin{equation}\label{eqn_faa}
f(g(x))= f(b_{0}) + \sum_{n=1}^{\infty}\sum_{k=1}^{n}f^{(k)}(b_{0})B_{n,k}(b_{1},\ldots,b_{n-k+1})\frac{x^{n}}{n!}.
\end{equation}

\section{Differential equation $y^{\prime}=f(y)$ and Bell polynomials}

This section contains the general results of this paper. Here we show the condition that the function $g(x)$ must satisfy for $g(t+g^{-1}(x))$ to be a solution of the differential equation $y^{\prime}=f(y)$. Then we will give a representation of $g(t+g^{-1}(x))$ in power series using partial Bell polynomials and finally we use the solution of the equation $y^{\prime}=f(y)$ to find special values of partial Bell polynomials.

\begin{theorem}
The function $y(t,x)=g(t+g^{-1}(x))$ is solution of the differential equation $y^{\prime}=f(y)$ with initial value problem $y(0)=x$, where $f(x)=\frac{1}{(g^{-1})^{\prime}(x)}$.
\end{theorem}
\begin{proof}
Using the method of separation of variables and the value of the function $f(x)$ given in the hypothesis, we find that
\begin{eqnarray*}
    \int\frac{dy}{f(y)}&=&\int(g^{-1})^{\prime}(y)dy\\
    &=&g^{-1}(y)=t+C.
\end{eqnarray*}
As $y(0)=x$, then $y(t,x)=g(t+g^{-1}(x))$.
\end{proof}

\begin{theorem}\label{theo_prin}
The function $g(t+g^{-1}(x))$ has the following representation using Bell polynomials
\begin{multline}
    g(t+g^{-1}(x))
    =x+g^{\prime}[g^{-1}(x)]t\\+\sum_{n=1}^{\infty}\sum_{k=1}^{n}f^{(k)}(x)B_{n,k}(g^{\prime}[g^{-1}(x)],\ldots,g^{(n-k+1)}[g^{-1}(x)])\frac{t^{n+1}}{(n+1)!},
\end{multline}
where $g^{(n)}[g^{-1}(x)]$ satisfies
\begin{equation}\label{eqn_coeff_n}
    g^{(n+1)}[g^{-1}(x)]=\sum_{k=1}^{n}f^{(k)}(x)B_{n,k}(g^{\prime}[g^{-1}(x)],\ldots,g^{(n-k+1)}[g^{-1}(x)]),\ \ \ n\geq1
\end{equation}
with $g^{-1}(y)=\int\frac{dy}{f(y)}$.
\end{theorem}
\begin{proof}
Applying Taylor formula to $g(t+g^{-1}(x))$ we get
\begin{eqnarray*}
    g(t+g^{-1}(x))&=&x+\sum_{n=1}^{\infty}g^{(n)}(g^{-1}(x))\frac{t^{n}}{n!}\\
    &=&x+g^{\prime}[g^{-1}(x)]t+\sum_{n=1}^{\infty}g^{(n+1)}(g^{-1}(x))\frac{t^{n+1}}{(n+1)!}.
\end{eqnarray*}
Since $g(t+g^{-1}(x))$ is solution of $y^{\prime}=f(y)$ with initial value problem $y(0)=x$, then by directly applying Fa\`a di Bruno's formula (\ref{eqn_faa}) to $y^{\prime}=f(y)$ we obtain the desired result.
\end{proof}

\begin{theorem}\label{theo_sec}
Let $y^{\prime}=f(y)$ be the autonomous differential equation with initial value problem $y(0)=x$. For $n\geq k\geq1$ we have
\begin{multline}
    B_{n,k}(g^{\prime}[g^{-1}(x)],g^{\prime\prime}[g^{-1}(x)],\ldots,g^{(n-k+1)}[g^{-1}(x)])\\
    =\frac{1}{k!}\sum_{i=0}^{k}\binom{k}{i}(-1)^{k-i}x^{k-i}\frac{d^{n}}{dt^{n}}g^{i}(t+g^{-1}(x))\bigg\vert_{t=0}.
\end{multline}
\end{theorem}
\begin{proof}
Making $x_{m}=g^{(m)}[g^{-1}(x)]$ in the equation (\ref{eqn_egf_pow}) leads us to
\begin{multline*}
    \sum_{n=k}^{\infty}B_{n,k}(g^{\prime}[g^{-1}(x)],\ldots,g^{(n-k+1)}[g^{-1}(x)])\frac{t^{n}}{n!}\\
    =\frac{1}{k!}\left(\sum_{m=1}^{\infty}g^{(m)}[g^{-1}(x)]\frac{t^{m}}{m!}\right)^{k}\\
    =\frac{1}{k!}(g(t+g^{-1}(x))-x)^{k}\\
    =\frac{1}{k!}\sum_{i=0}^{k}\binom{k}{i}(-1)^{k-i}x^{k-i}g^{i}(t+g^{-1}(x)).
\end{multline*}
Differentiating $m\geq k\geq0$ times with respect to $t$
\begin{multline*}
    \sum_{n=k}^{\infty}B_{n,k}(g^{\prime}[g^{-1}(x)],\ldots,g^{(n-k+1)}[g^{-1}(x)])\frac{t^{n-m}}{(n-m)!}\\
    =\frac{1}{k!}\sum_{i=0}^{k}\binom{k}{i}(-1)^{k-i}x^{k-i}\frac{d^{m}}{dt^{m}}g^{i}(t+g^{-1}(x))
\end{multline*}
and then by making $t=0$, we obtain the desired result
\[
B_{m,k}(g^{\prime}[g^{-1}(x)],\ldots,g^{(m-k+1)}[g^{-1}(x)])
    =\frac{1}{k!}\sum_{i=0}^{k}\binom{k}{i}(-1)^{k-i}x^{k-i}\frac{d^{m}}{dt^{m}}g^{i}(t+g^{-1}(x))\bigg\vert_{t=0}.
    \qedhere
\]
\end{proof}

\section{Some examples}

We will use the Theorem \ref{theo_prin} to represent the solution $g(t+g^{-1}(x))$ by Bell polynomials when $f(x)$ is any of the following functions: $f(x)=ax$, $f(x)=e^{x}$, $f(x)=\sqrt{1-x^2}$, $f(x)=1+x^{2}$, $f(x)=\sqrt{x^{2}\pm1}$, $f(x)=1-x^2$ and $f(x)=\frac{1}{(x-1)^{\alpha-1}}$, where $a,\alpha\in\C$. In addition by using the Theorem \ref{theo_sec} we will find identities for Bell polynomials, some of which were constructed by Feng Qi et al in \cite{FengQi1}, \cite{FengQi2}, \cite{FengQi3}. In particular, we will note that we can associate Stirling numbers and Lah numbers with autonomous differential equations of order one.

\subsection{Equation $y^{\prime}=ay$}

The solution to this equation is $y(t,x)=xe^{at}$ where $f(x)=ax$, $f^{(k)}(x)=0$ for $k\geq2$, and $g^{(n)}[g^{-1}(x)]=a^{n}x$. By the Theorem \ref{theo_prin}, we have the representation of $xe^{at}$ using Bell polynomials
\begin{eqnarray*}
    xe^{ax}&=&x+axt+a\sum_{n=1}^{\infty}\sum_{k=1}^{n}f^{(k)}(x)B_{n,k}(ax,a^{2}x,\ldots,a^{n-k+1}x)\frac{t^{n+1}}{(n+1)!}\\
    &=&x+axt+a\sum_{n=1}^{\infty}B_{n,1}(ax,a^2x,\ldots,a^{n}x)\frac{t^{n+1}}{(n+1)!}
\end{eqnarray*}
and by the Theorem \ref{theo_sec} we find the following value of Bell polynomial using the solution $y(t,x)=xe^{at}$
\begin{eqnarray*}
    B_{n,k}(ax,a^{2}x,\ldots,a^{n-k+1}x)&=&\frac{1}{k!}\sum_{i=0}^{k}\binom{k}{i}(-1)^{k-i}x^{k-i}\frac{d^n}{dt^n}x^{i}e^{ait}\bigg\vert_{t=0}\\
    &=&\frac{a^{n}x^{k}}{k!}\sum_{i=0}^{k}\binom{k}{i}(-1)^{k-i}i^{n}\\
    &=&a^{n}x^{k}\stirlingsk{n}{k},
\end{eqnarray*}
where we have used that
\begin{equation*}
    \stirlingsk{n}{k}=\frac{1}{k!}\sum_{i=0}^{k}\binom{k}{i}(-1)^{k-i}i^{n}.
\end{equation*}

We then relate the Stirling numbers of the second kind to the differential equation $y^{\prime}=ay$. When we make $x=1$, $a=1$, the above provides another proof for the known result $B_{n,k}(1,1,\ldots,1)=\stirlingsk{n}{k}$.

\subsection{Equation $y^{\prime}=e^{y}$}

The solution to this equation is $y(t,x)=-\ln(-t+e^{-x})$. Then by the Theorem \ref{theo_prin}, with $f(x)=e^{x}$ and $g^{(n)}[g^{-1}(x)]=(n-1)!e^{nx}$, we have the representation of the function $-\ln(-t+e^{-x})$, that is,
\begin{multline*}
    -\ln(-t+e^{-x})=x+e^{x}t\\
    +e^{x}\sum_{n=1}^{\infty}\sum_{k=1}^{n} B_{n,k}\left(0!e^{x},1!e^{2x},\ldots,(n-k)!e^{(n-k+1)x}\right)\frac{t^{n+1}}{(n+1)!}\\
    =x+e^{x}t
    +e^{x}\sum_{n=1}^{\infty}e^{nx}\sum_{k=1}^{n}B_{n,k}\left(0!,1!,\ldots,(n-k)!\right)\frac{t^{n+1}}{(n+1)!}\\
    =x+e^{x}t
    +e^{x}\sum_{n=1}^{\infty}e^{nx}B_{n}\left(0!,1!,\ldots,(n-1)!\right)\frac{t^{n+1}}{(n+1)!}.
\end{multline*}
By the Theorem \ref{theo_sec} we find another special value of Bell polynomials, i.e.
\begin{multline*}
    \sum_{n=k}^{\infty}B_{n,k}\left(0!e^{x},1!e^{2x},\ldots,(n-k+1)!e^{(n-k+1)x}\right)\frac{t^n}{n!}\\
    =\frac{(-\ln(-t+e^{-x})-x)^k}{k!}
    =\frac{(-\ln(-te^{x}+1))^k}{k!}
    =\sum_{n=k}^{\infty}e^{nx}\stirlingfk{n}{k}\frac{t^n}{n!}
\end{multline*}
and by comparing the coefficients of the two sums we get
\begin{equation*}
    B_{n,k}\left(0!e^{x},1!e^{2x},\ldots,(n-k+1)!e^{(n-k+1)x}\right)=e^{nx}\stirlingfk{n}{k}.
\end{equation*}
Clearly the unsigned Stirling numbers of the first kind are related to the differential equation $y^{\prime}=e^{y}$.

\subsection{Equation $y^{\prime}=\sqrt{1-y^{2}}$}

The solution of the equation with initial value problem $y(0)=x$ is
\begin{eqnarray*}
    y(t,x)&=&\sin(t+\arcsin(x))\\
    &=&x\cos(t)+\sqrt{1-x^2}\sin(t).
\end{eqnarray*}
Then, due to the Theorem \ref{theo_prin}, we reach
\begin{multline*}
    x\cos(t)+\sqrt{1-x^2}\sin(t)=x+\sqrt{1-x^2}t\\
    +\sum_{n=1}^{\infty}\sum_{k=1}^{n}f^{(k)}(x)B_{n,k}\left(\sqrt{1-x^2},-x,\ldots,-\cos\left[\arcsin(x)+(n-k)\frac{\pi}{2}\right]\right)\frac{t^{n+1}}{(n+1)!},
\end{multline*}
where
\begin{equation*}
    f^{(k)}(x)=\sum_{i=1}^{k}(-1)^{i}2^{2i-k}k!\binom{1/2}{i}\binom{i}{k-i}(1-x^2)^{1/2-i}x^{2i-k}
\end{equation*}
has been obtained by applying Fa\`a di Bruno's formula to $f(x)=g(h(x))$ with $g(x)=\sqrt{x}$ and $h(x)=1-x^2$ and by the equation (\ref{eqn_ax_a})
\begin{equation*}
    B_{n,k}(-2x,-x,0,\ldots,0)=(-1)^{k}2^{2k-n}\frac{n!}{k!}\binom{k}{n-k}x^{2k-n}.
\end{equation*}

A simple application of the previous result with $x=0$, $x=1$, $x=\frac{1}{2}$, $x=\frac{\sqrt{3}}{2}$ and $x=\frac{\sqrt{2}}{2}$ leads us to the representation of the functions
$\sin(t)$, $\cos(t)$, $\frac{1}{2}\cos(t)+\frac{\sqrt{3}}{2}\sin(t)$, $\frac{\sqrt{3}}{2}\cos(t)+\frac{1}{2}\sin(t)$ and $\frac{\sqrt{2}}{2}(\cos(t)+\sin(t))$, respectively.

Feng Qi \cite{FengQi2} showed the following result
\begin{multline}\label{eqn_bell_diff_sin}
    B_{n,k}(\cos x,-\sin x,-\cos x,\sin x,\ldots,-\cos\left[x+(n-k)\frac{\pi}{2}\right])\\
    =\frac{(-1)^{k}}{k!}\sum_{i=0}^{k}\frac{1}{2^{i}}\binom{k}{i}\sin^{k-i}x\sum_{j=0}^{i}(-1)^{j}\binom{i}{j}(2j-i)^{n}\\
    \times\cos\left[(n-i)\frac{\pi}{2}+(2j-i)x\right],
\end{multline}
which when composing with $\arcsin x$ leads to the following result
\begin{multline}
    B_{n,k}\left(\sqrt{1-x^2},-x,-\sqrt{1-x^2},x,\ldots,-\cos\left[\arcsin(x)+(n-k)\frac{\pi}{2}\right]\right)\\
    =\frac{(-1)^{k}}{k!}\sum_{i=0}^{k}\binom{k}{i}\frac{x^{k-i}}{2^{i}}\sum_{j=0}^{i}(-1)^{j}\binom{i}{j}(2j-i)^{m}\\
    \times\cos\left[(n-i)\frac{\pi}{2}+(2j-i)\arcsin(x)\right].
\end{multline}

Analogously, we can find similar results for the solution $\cos(t+\arccos(x))$ of the differential equation $y^{\prime}=-\sqrt{1-y^2}$, with $y(0)=x$.

\subsection{Equation $y^{\prime}=1+y^{2}$}

The solution of this equation is $\tan(t+\arctan(x))$ with $f(x)=1+x^2$. By the Theorem~\ref{theo_sec} it follows that
\begin{theorem}
\begin{multline}
    \tan^{(n+1)}(\arctan(x))=2x\tan^{(n)}(\arctan(x))\\
    +\sum_{k=1}^{n-1}\binom{n}{k}\tan^{(k)}(\arctan(x))\tan^{(n-k)}(\arctan(x))
\end{multline}
\end{theorem}
\begin{proof}
It is clear that $f^{(k)}(x)=0$ for all $k\geq3$. Then the theorem follows by keeping in mind that
\begin{equation*}
    B_{n,1}(\tan^{\prime}(\arctan(x)),\ldots,\tan^{(n)}(\arctan(x)))=\tan^{(n)}(\arctan(x))
\end{equation*}
and
\begin{align*}
    B_{n,2}&(\tan^{\prime}(\arctan(x)),\ldots,\tan^{(n-1)}(\arctan(x))) \\
    &{ }=\frac{1}{2}\sum_{k=1}^{n-1}\binom{n}{k}\tan^{(k)}(\arctan(x))\tan^{(n-k)}(\arctan(x)).
\qedhere
\end{align*}
\end{proof}

Now using the Theorem \ref{theo_prin} together with the previous result we obtain
\begin{theorem}
The representation of the function $\tan(t+\arctan(x))$ is
\begin{multline*}
    \tan(t+\arctan(x))=x+(1+x^2)t+\\
    =x+(1+x^2)t+2x(1+x^2)\frac{t^2}{2!}+2x\sum_{n=2}^{\infty}\tan^{(n)}(\arctan(x))\frac{t^{n+1}}{(n+1)!}\\
    +\sum_{n=2}^{\infty}\sum_{k=1}^{n-1}\binom{n}{k}\tan^{(k)}(\arctan(x))\tan^{(n-k)}(\arctan(x))\frac{t^{n+1}}{(n+1)!},
\end{multline*}
where composing the equation (1.12) in \cite{FengQi2} with $\arctan(x)$ leads us to
\begin{multline*}
    \tan^{(n)}(\arctan(x))=
    -\sum_{k=1}^{n+1}\frac{1}{k}\sum_{l=0}^{k}\frac{(-1)^{l}}{2^{l}}\binom{k}{l}(x^{2}+1)^{l/2}\\
    \times\sum_{q=0}^{l}\binom{l}{q}(2q-l)^{n+1}\sin\left[\frac{\pi}{2}n+(2q-l)\arctan(x)\right].
\end{multline*}
\end{theorem}

Analogous results are obtained for the solution $\cot(t+\mathrm{arccot}(x))$ of the differential equation $y^{\prime}=-1-y^2$ with $y(0)=x$.

\subsection{Equation $y^{\prime}=\sqrt{y^{2}\pm1}$}

The solution of the equation $y^{\prime}=\sqrt{y^{2}+1}$ with initial value problem $y(0)=x$ is
\begin{eqnarray*}
    y(t,x)&=&\sinh(t+\mathrm{arcsinh}(x))\\
    &=&x\cosh(t)+\sqrt{1+x^2}\sinh(t).
\end{eqnarray*}
For $n\geq1$ denote $\phi_{n}(x)$ the function
\begin{equation}
    \phi_{n}(x)=
    \begin{cases}
        x,& \text{ if $n$ is an even number;}\\
        \sqrt{1+x^2}, & \text{ if $n$ is an odd number.}
    \end{cases}
\end{equation}

By the Theorem \ref{theo_prin} we have
\begin{multline*}
    x\cosh(t)+\sqrt{1+x^2}\sinh(t)=x+\sqrt{1+x^2}t\\
    +\sum_{n=1}^{\infty}\sum_{k=1}^{n}f^{(k)}(x)B_{n,k}\left(\phi_{1}(x),\phi_{2}(x),\ldots,\phi_{n-k+1}(x)\right)\frac{t^{n+1}}{(n+1)!},
\end{multline*}
where
\begin{equation*}
    f^{(k)}(x)=\sum_{i=1}^{k}2^{2i-k}k!\binom{1/2}{i}\binom{i}{k-i}(1+x^2)^{1/2-i}x^{2i-k}.
\end{equation*}
has been obtained by applying Fa\`a di Bruno's formula to $f(x)=g(h(x))$ with $g(x)=\sqrt{x}$ and $h(x)=1+x^2$ and by the equation (\ref{eqn_ax_a})
\begin{equation*}
    B_{n,k}(2x,2,0,\ldots,0)=2^{2k-n}\frac{n!}{k!}\binom{k}{n-k}x^{2k-n}.
\end{equation*}

\begin{theorem}
For $n\geq k\geq1$ we have
\begin{multline}\label{eqn_sol_sinh}
    B_{n,k}(\phi_{1}(x),\phi_{2}(x),\ldots,\phi_{n-k+1}(x))\\
    =\frac{(-1)^{k}x^k}{k!}\sum_{i=0}^{k}\binom{k}{i}\frac{1}{2^i}\sum_{l=0}^{i}\binom{i}{l}(-1)^{l}(2l-i)^{n}\left(\sqrt{1+x^2}+x\right)^{2l-i}.
\end{multline}
\end{theorem}
\begin{proof}
Using $\sinh(t)=\frac{e^{t}-e^{-t}}{2}$ we will obtain $\frac{d^{n}}{dt^n}\sinh^i(t)$, that is,
\begin{eqnarray*}
\frac{d^{n}\sinh^i(t)}{dt^n}&=&\frac{d^{n}}{dt^n}\frac{(e^{t}-e^{-t})^i}{2^i}\\
&=&\frac{d^{n}}{dt^n}\frac{1}{2^i}\sum_{l=0}^{i}\binom{i}{l}(-1)^{i-l}e^{(2l-i)t}\\
&=&\frac{1}{2^i}\sum_{l=0}^{i}\binom{i}{l}(-1)^{i-l}(2l-i)^{n}e^{(2l-i)t}.
\end{eqnarray*}
Then by the Theorem \ref{theo_sec} 
\begin{multline*}
    B_{n,k}(\phi_{1}(x),\phi_{2}(x),\ldots,\phi_{n-k+1}(x))\\
    =\frac{1}{k!}\sum_{i=0}^{k}\binom{k}{i}(-1)^{k-i}x^{k-i}\frac{1}{2^i}\sum_{l=0}^{i}\binom{i}{l}(-1)^{i-l}(2l-i)^{n}e^{(2l-i)\mathrm{arcsinh}(x)}\\
    =\frac{(-1)^{k}x^k}{k!}\sum_{i=0}^{k}\binom{k}{i}\frac{1}{2^i}\sum_{l=0}^{i}\binom{i}{l}(-1)^{l}(2l-i)^{n}\left(\sqrt{1+x^2}+x\right)^{2l-i}
\end{multline*}
and the theorem follows.
\end{proof}

Analogously, we can find similar results to the previous ones for the solution $\cosh(t+\mathrm{arccosh}(x))$ of the differential equation $y^{\prime}=\sqrt{y^2-1}$, with $y(0)=x$, changing $x^{2}+1$ into $x^{2}-1$ and $\phi_{n}(x)$ into $\psi_{n}(x)$, where
\begin{equation}
    \psi_{n}(x)=
    \begin{cases}
        x,&\text{ if $n$ is even; }\\
        \sqrt{x^{2}-1},&\text{ if $n$ is odd}.
    \end{cases}
\end{equation}

\subsection{Equation $y^{\prime}=1-y^{2}$}

The solution of this equation is $\tanh(t+\mathrm{arctanh}(x))$. Before applying the Theorems~\ref{theo_prin} and \ref{theo_sec} to the differential equation, we first calculate the $n$-th derivative of the function $\tanh(x)$.
\begin{theorem}\label{theo_der_tanh}
\begin{equation*}
    \frac{d^n\tanh(x)}{dt^n}=
    -\sum_{k=1}^{n+1}   \frac{\tanh^k(x)}{k}\sum_{i=0}^{k}\binom{k}{i}\frac{1}{2^i}\sum_{l=0}^{i}\binom{i}{l}(-1)^{l}(2l-i)^{n}e^{(2l-i)x}.
\end{equation*}
\end{theorem}
\begin{proof}
\begin{multline*}
\frac{d^n\tanh(x)}{dt^n}=\frac{d^{n+1}\ln(\cosh(x))}{dt^{n+1}}\\
    =\sum_{k=1}^{n+1}(\ln u)^{(k)}B_{n+1,k}(\psi_{1}(\sinh(x)),\psi_{2}(\sinh(x)),\ldots,\psi_{n-k+2}(\sinh(x)))\\
    =\sum_{k=1}^{n+1}\frac{(-1)^{k-1}(k-1)!}{\cosh^k(x)} B_{n+1,k}(\psi_{1}(\sinh(x)),\psi_{2}(\sinh(x)),\ldots,\psi_{n-k+2}(\sinh(x))).
\end{multline*}
Changing $x^{2}+1$ into $x^{2}-1$ and $\phi_{n}(x)$ into $\psi_{n}(x)$ in the equation (\ref{eqn_sol_sinh}) and then composing with $\sinh(x)$ we reach
\begin{equation*}
    \frac{d^n\tanh(x)}{dt^n}=
    -\sum_{k=1}^{n+1}   \frac{\tanh^k(x)}{k}\sum_{i=0}^{k}\binom{k}{i}\frac{1}{2^i}\sum_{l=0}^{i}\binom{i}{l}(-1)^{l}(2l-i)^{n}e^{(2l-i)x}.
\qedhere
\end{equation*}
\end{proof}

By the Theorem \ref{theo_prin} we obtain the following result
\begin{theorem}
\begin{multline}
    \tanh^{(n+1)}(\mathrm{arctanh}(x))=-2x\tanh^{(n)}(\mathrm{arctanh}(x))\\
    -\sum_{k=1}^{n-1}\binom{n}{k}\tanh^{(k)}(\mathrm{arctanh}(x))\tanh^{(n-k)}(\mathrm{arctanh}(x)).
\end{multline}
\end{theorem}

Finally, by the Theorems \ref{theo_prin} and \ref{theo_der_tanh}
\begin{theorem}
The representation of the function $\tanh(t+\mathrm{arctanh}(x))$ is
\begin{multline*}
    \tanh(t+\mathrm{arctanh}(x))\\
    =x+(1-x^2)t-2x(1-x^2)\frac{t^2}{2!}-2x\sum_{n=2}^{\infty}\tanh^{(n)}(\mathrm{arctanh}(x))\frac{t^{n+1}}{(n+1)!}\\
    -\sum_{n=2}^{\infty}\sum_{k=1}^{n-1}\binom{n}{k}\tanh^{(k)}(\mathrm{arctanh}(x))\tanh^{(n-k)}(\mathrm{arctanh}(x))\frac{t^{n+1}}{(n+1)!},
\end{multline*}
where
\begin{equation*}
    \tanh^{(n)}(\mathrm{arctanh}(x))=-\sum_{k=1}^{n+1}   \frac{x^k}{k}\sum_{i=0}^{k}\binom{k}{i}\frac{1}{2^i}\sum_{l=0}^{i}\binom{i}{l}(-1)^{l}(2l-i)^{n}\left(\frac{\sqrt{1+x}}{\sqrt{1-x}}\right)^{2l-i}.
\end{equation*}
\end{theorem}

\subsection{Equation $y^{\prime}=\frac{1}{(y\pm1)^{\alpha-1}}$, $\alpha\in\C$}

The solution of the equation is
\begin{equation}\label{eqn_diff_1}
    y(t,x)=
    \begin{cases}
        \mp1+(x\pm1)\left(1+\frac{\alpha t}{(x\pm1)^{\alpha}}\right)^{1/\alpha},&\ \text{if $\alpha\neq0$,}\\
        \mp1+(x\pm1)e^{t},&\ \text{if $\alpha=0$}.
    \end{cases}
\end{equation}

Assume $\alpha\neq0$, since the case $\alpha=0$ is obtained from the solution of $y^{\prime}=ay$. Then by the Theorem \ref{theo_prin} we get the representation
\begin{multline}
\mp1+(x\pm1)\left(1+\frac{\alpha t}{(x\pm1)^{\alpha}}\right)^{1/\alpha}=\\
x+\frac{t}{(x\pm1)^{\alpha-1}}+
    \sum_{n=1}^{\infty}\sum_{k=1}^{n}\frac{(-1)^{k}\ffact{\alpha+k-2}_{k}}{(x\pm1)^{\alpha +k-1}}\\
    \times B_{n,k}\left(\frac{\alpha}{(x\pm1)^{\alpha-1}}\ffactb{\frac{1}{\alpha}}_{1},\ldots,\frac{\alpha^{n-k+1}}{(x\pm1)^{\alpha(n-k+1)-1}}\ffactb{\frac{1}{\alpha}}_{n-k+1}\right)\frac{t^{n+1}}{(n+1)!},
\end{multline}
where
\begin{equation}
    \ffact{a}_{n}=\prod_{k=0}^{n-1}(a-k)=
    \begin{cases}
        a(a-1)\cdots(a-n+1),&\text{ if $n\geq1$; }\\
        1,&\text{ if $n=0$,}
    \end{cases}
\end{equation}
is called the falling factorial.

Now we apply (\ref{eqn_hom}) to obtain
\begin{multline}\label{eqn_bell_fal}
    B_{n,k}\left(\frac{\alpha}{(x\pm1)^{\alpha-1}}\ffactb{\frac{1}{\alpha}}_{1},\ldots,\frac{\alpha^{n-k+1}}{(x\pm1)^{\alpha(n-k+1)-1}}\ffactb{\frac{1}{\alpha}}_{n-k+1}\right)\\
    =\frac{\alpha^n}{(x\pm1)^{\alpha n-k}}B_{n,k}\left(\ffactb{\frac{1}{\alpha}}_{1},\ldots,\ffactb{\frac{1}{\alpha}}_{n-k+1}\right).
\end{multline}
Then by Theorem 2.1 in \cite{FengQi1}
\begin{multline}\label{eqn_bell_fal_fin}
    B_{n,k}\left(\frac{\alpha}{(x\pm1)^{\alpha-1}}\ffactb{\frac{1}{\alpha}}_{1},\ldots,\frac{\alpha^{n-k+1}}{(x\pm1)^{\alpha(n-k+1)-1}}\ffactb{\frac{1}{\alpha}}_{n-k+1}\right)\\
    =\frac{(-1)^{k}}{k!}\frac{\alpha^{n}}{(x\pm1)^{\alpha n-k}}\sum_{i=0}^{k}(-1)^{i}\binom{k}{i}\ffactb{\frac{i}{\alpha}}_{n}.
\end{multline}
Now by the equation (\ref{eqn_coeff_n})
\begin{theorem}
\begin{equation}\label{eqn_iden_ffact}
    \sum_{k=1}^{n}\frac{\ffact{\alpha+k-2}_{k}}{k!}\sum_{i=0}^{k}(-1)^{i}\binom{k}{i}\ffactb{\frac{i}{\alpha}}_{n}=\alpha\ffactb{\frac{1}{\alpha}}_{n+1}.
\end{equation}
\end{theorem}

Finally we will note the relationship that exists between the solution of the differential equation $y^{\prime}=(y+1)^2$, with $y(0)=0$, and the Lah numbers. We will make $\alpha=-1$ in (\ref{eqn_bell_fal}) to obtain
\begin{multline*}
    B_{n,k}(-\ffact{-1}_{1},\ffact{-1}_{2},\ldots,(-1)^{n-k+1}\ffact{-1}_{n-k+1})\\
    =(-1)^{n}B_{n,k}(-1!,2!,-3!,\ldots,(-1)^{n-k+1}(n-k+1)!)\\
    =B_{n,k}(1!,2!,3!,\ldots,(n-k+1)!)=\binom{n-1}{k-1}\frac{n!}{k!}.
\end{multline*}
Then by (\ref{eqn_bell_fal_fin}) it is proved that
\begin{equation}
    \binom{n-1}{k-1}\frac{n!}{k!}=\frac{(-1)^{n+k}}{k!}\sum_{i=0}^{k}(-1)^{i}\binom{k}{i}\ffact{-i}_{n}
\end{equation}
and in combination with the equation (\ref{eqn_iden_ffact})
\begin{theorem}
For $n\geq1$
\begin{equation}
    \sum_{k=1}^{n}(-1)^{k}\ffact{k-3}_{k}\binom{n-1}{k-1}\frac{n!}{k!}=(n+1)!.
\end{equation}
\end{theorem}

\section{Conclusion}
As noted, Fa\`a di Bruno's formula allows us to connect Bell polynomials with autonomous differential equations of order one. In particular, we find differential equations for the Stirling and Lah numbers. The same procedure should allow us to find differential equations for the vast amount of Bell polynomial values found in the existing literature.

\subsection*{Acknowledgment}
The author is thankful to the anonymous referee for the helpful comments and suggestions which helped improve the paper.


\EditInfo{June 30, 2021}{September 06, 2021}{Karl Dilcher}

\end{paper}